\documentclass{amsart}

\usepackage{latexsym,amsmath,amsfonts,amscd,amssymb, amsthm}

\theoremstyle{plain}
\newtheorem{theorem}{Theorem}[section]

\theoremstyle{definition}
\newtheorem{definition}{Definition}[section]
\newtheorem{example}{Example}[section]

\theoremstyle{remark}
\newtheorem{remark}{Remark}[section]

\title{A remark on finite type conditions}

\author{John P. D'Angelo}

\address{Dept. of Mathematics, Univ. of Illinois, 1409 W. Green St., Urbana IL 61801}

\email{jpda@math.uiuc.edu}

\begin{document}

\maketitle

\begin{abstract} We prove that a certain positivity condition, considerably more general
than pseudoconvexity, enables one to conclude that the regular order of contact and singular order
of contact agree when these numbers are $4$.

\medskip

\noindent
{\bf AMS Classification Numbers}: 32F18, 32T25, 32T27, 32U05, 32V35.

\medskip

\end{abstract}

\section{introduction}

Let $(M,p)$ be the germ of a smooth real submanifold of complex Euclidean space ${\mathbb C}^n$. 
It is natural to ask whether there are any positive dimensional complex analytic germs $(V,p)$ contained in $(M,p)$,
and if not, how close such varieties can contact $M$ at $p$. When $M$ is a real hypersurface,
the book [D] discusses this problem and its relationship to estimates for the ${\overline \partial}$-Neumann problem.
See also [C], [D1], [DF], and [K]. This paper returns to these matters and corrects an error in [D],
where the hypothesis of pseudoconvexity was omitted in the statement of a minor result.

One of the crucial issues in the above discussion involves singularities.
Suppose for example that $(V,0)$ is the germ at the origin of an irreducible
complex analytic $1$-dimensional variety in ${\mathbb C}^n$
that is singular at $0$. Assume $V$ is locally defined by holomorphic functions
$f_1,...,f_k$ and consider the real hypersurface $M$ in ${\mathbb C}^{n+1}$ defined near $0$ by the equation
$$ 2{\rm Re}(z_{n+1}) + \sum_{j=1}^k |f_j(z)|^2 = 0. $$
Then $M$ contains the complex variety defined by 
$$ z_{n+1}=f_1(z)=...= f_k(z) = 0. $$
The germ $(M,0)$ thus contains the singular holomorphic curve $(V,0)$, but by the above assumptions, no non-singular one.
The general theory from [D] and [D1] therefore must consider singularities.

One naturally asks under what geometric information 
it suffices to consider only non-singular curves. 
McNeal [M] and Boas-Straube [BS] showed, when $M$ bounds a {\it convex} domain, 
that one need consider only smooth complex varieties in the discussion.
The same conclusion holds for pseudoconvex Reinhardt domains ([FIK]).

Let $\Delta_q(M,p)$ denote the maximum order of contact of a $q$-dimensional complex analytic
variety with $M$ at $p$, as defined in [D1]. Let $\Delta_q^{\rm reg} (M,p)$ denote the maximum 
order of contact of a $q$-dimensional complex analytic manifold with $M$ at $p$. When $M$ bounds a domain
that is convex near $p$ or bounds a pseudoconvex Reinhardt domain, the above authors showed that 
$$\Delta_q (M,p) = \Delta_q^{\rm reg} (M,p). $$ 
Kohn [K] had noted this equality when $M$ is pseudoconvex and $\Delta_q^{\rm reg} (M,p) = 2$.

\begin{remark} There are several distinct ways to define the order of contact of a singular complex variety
of dimension larger than $1$  with a real hypersurface. When $q>1$, the number $\Delta_q(M,p)$ need
not equal the measurement defined in [C], although the measurements are simultaneously finite.
See [BN], [BN1],  and [F] for examples and inequalities relating these numbers. 
 \end{remark}

In this paper we consider only the case when $q=1$. In [D] the author stated that
$\Delta_1 (M,p) = \Delta_1^{\rm reg} (M,p) $ when $\Delta_1^{\rm reg} (M,p) = 4$. The author had intended
to assume $M$ was pseudoconvex near $p$, but incorrectly omitted this hypothesis. McNeal and Mernik
[MM] gave an example of a hypersurface $M$ in ${\mathbb C}^3$, defined by a polynomial equation of degree $5$, where
$\Delta_1^{\rm reg} (M,0) = 4$, but $M$ contains a singular complex analytic curve through $0$, and hence
$\Delta_1 (M,0) = \infty$. They also proved, when $M$ is pseudoconvex near $0$, that
$\Delta_1^{\rm reg} (M,0) = 4$ implies $\Delta_1 (M,0) = 4$. 

In order to state the main result of this paper, we recall a positivity property {\bf PS} for functions from [D1].
We then extend this concept to a positivity property,
considerably weaker than pseudoconvexity, for germs of real hypersurfaces.
We provide a simple proof that this more general property implies the above equality.

We use the words {\it pure terms} for any harmonic polynomial and {\it mixed terms} for any sum
of monomials that are neither holomorphic nor anti-holomorphic.

\begin{definition} Let $g:({\mathbb C}^n,0) \to ({\mathbb R}, 0)$ be the germ of a smooth function.
We say that $g$ satisfies property {\bf PS} if, whenever $z:({\mathbb C},0) \to ({\mathbb C}^n,0)$
is the germ of a holomorphic map for which the pullback $z^*g$ vanishes to finite order, the order $2k$ is even and
the Taylor coefficient of $|t|^{2k}$ is positive. \end{definition}

In Definition 1.1 we allow the function $z^*g$ to vanish to infinite order. In this paper, however, we will always
be in the {\it finite type} situation; every pullback to a non-constant map vanishes to finite order.

Property {\bf PS} has a simple interpretation in terms of the Laplacian.
Let $t$ be a complex variable. Put $t=|t|e^{i\theta}$.
The operator $L= {d \over dt} \ {d \over d{\overline t}}$
is a constant times the usual Laplacian. Let $u$ be a smooth function defined near $0$ in ${\mathbb C}$ with $u(0)=0$.
If $u$ vanishes to even order at $0$, then the lowest order terms in its Taylor series can be written
$$ \sum_{j=0}^{2k} c_j t^j {\overline t}^{2k-j} = |t|^{2k} p(\theta). \eqno (1) $$
Here $p$ is a trig polynomial.
Property {\bf PS} guarantees that $(L^k u)(0) = (k!)^2 c_k > 0$. The coefficient $c_k$ is the average value
of $p$ on the circle.  A strictly subharmonic function which vanishes to order two at $0$
satisfies {\bf PS} with $k=1$. A smooth subharmonic function that vanishes
to finite order $2k$ at $0$, and whose Taylor expansion up to that order has no pure terms, satisfies {\bf PS} as well.

If $r$ is plurisubharmonic, and $z$ is as in Definition 1.1, then $z^*r$ is subharmonic.
Thus plurisubharmonic functions with no pure terms satisfy property {\bf PS}.
This property depends only on the Taylor series of $r$ at $0$, whereas plurisubharmonicity
depends on the values of $r$ in a neighborhood of $0$.

\begin{remark} We use the term ``property {\bf PS}'' for the following reason.
In [D1] the author called essentially the same concept ``property P''.
The concept was defined there only for functions, rather than for hypersurfaces.
Later Catlin used the term ``property P'' for a completely different notion
that implies compactness for the ${\overline \partial}$-Neumann problem.
The term, as used by Catlin, has become standard in subsequent work of many authors.
See [BN] for a brief mention of this matter.  \end{remark}

We write $C^\infty_p$ for the ring of germs of smooth functions at $p$. 
Let $(M,0)$ be the germ of a smooth real hypersurface in ${\mathbb C}^n$ and
let $r$ be a generator of the principal ideal in $C^\infty_0$ of functions
vanishing on $M$. We refer to $r$ as a {\it defining function} for $(M,0)$.
Let ${\bf j}_k r$ denote the $k$-th order Taylor polynomial of $r$ at the origin.
We wish to define property {\bf PS} for the germ $(M,0)$.

Since $(M,0)$ is a hypersurface, $dr(0)\ne 0$.
For $k\ge 2$ we can write 
$$ {\bf j}_k r = 2{\rm Re}(h_k) + g_k \eqno (2) $$
where $h_k$ is a holomorphic polynomial with $dh_k(0) \ne 0$
and $g_k$ is a polynomial containing only mixed terms.
We want the restriction of $g_k$ to the complex hypersurface defined by $h_k=0$ to satisfy {\bf PS}.
We denote this restriction by $G_k$.
When $(M,0)$ is real-analytic, we can write 
$$ r = 2{\rm Re}(h) + g $$
with $h$ a holomorphic germ and $g$ a real-analytic germ. In the $C^\infty$ case, however, $h$ and $g$
become formal power series. We avoid this problem as follows.

\begin{definition} Let $(M,0)$ be the germ of a smooth real hypersurface in ${\mathbb C}^n$.
We say that $(M,0)$ satisfies {\bf PS} if the following holds.
There is an integer $k_0$ such that, whenever $k \ge k_0$,
we can find a defining function $r_k$ such that (2) holds,
and the function $G_k$ satisfies ${\bf PS}$. 
\end{definition}

The reader might wonder about the meaning of the stabilization
condition in Definition 1.2. The next example and subsequent remark illustrate the idea.

\begin{example} Consider the function $r$ given for $m \ge 3$ by
$$  r(z) = 2 {\rm Re}(z_3) + |z_1 + z_2^m|^2. $$
After setting $z_3=0$, the restriction satisfies {\bf PS} at $0$, but its $k$-th order Taylor polynomial
does not when $k=2m-1$. We need $k \ge 2m$ to ensure that {\bf PS} holds.
\end{example}

\begin{remark} Let $(M_k,0)$ denote the germ defined by (2) and let $J_k$ denote the ideal
in $C^\infty_0$ defined by ${\bf j}_k r$. Recall from [D] that $(M,0)$ is {\it finite type} if and only if there is
an integer $k_0$ such that $\Delta_1(J_k) = \Delta_1(J_{k_0})$ for $k \ge k_0$, where $\Delta$
is defined in Definition 2.1 below. Thus {\it finite type} is a {\it finitely determined} condition.
In this paper we are concerned only with the simple case of hypersurface germs of type $4$.
In this case Theorem 1.1 enables us to ignore singular curves.
\end{remark}

We give one more class of functions whose germs satisfy {\bf PS}.

\begin{example} Put $r(z) = 2{\rm Re}(z_n) + g(\zeta)$ where $z=(\zeta,z_n)$
and $g$ is $C^\infty$. Assume its Taylor series contains no pure terms. 
Then $(M,0)$ satisfies {\bf PS} whenever $g$ does, 
such as when $g$ is plurisubharmonic. For example, $(M,0)$ satisfies {\bf PS}
whenever
$$ r(z) = 2{\rm Re}(z_n) = \sum_{j=1}^K |f_j(z)|^2 $$
for germs of holomorphic functions $f_j$ with $f_j(0)=0$.
\end{example}

\begin{remark} By Proposition 2 on Page 138 of [D], when $M$ is pseudoconvex near $0$, property
{\bf PS} holds for $(M,0)$. The conclusion of 
Theorem 1.1 therefore follows when $M$ is pseudoconvex near $0$. \end{remark}

\begin{theorem} Let $(M,0)$ be the germ of a smooth real hypersurface in ${\mathbb C}^n$
satisfying property {\bf PS}. If $\Delta_1^{\rm reg} (M,p)= 4$, then
$\Delta_1 (M,p) = 4$. \end{theorem}

\begin{remark} It is easy to show, without assuming {\bf PS}, that
$\Delta_1^{\rm reg} (M,p)= 3$ implies $\Delta_1(M,p)= 3$. It is well-known
that $\Delta_1^{\rm reg} (M,p)= 2$ implies $\Delta_1(M,p)= 2$.
See for example [K]. No such conclusion is possible when $\Delta_1^{\rm reg} (M,p)\ge 6$.
For example, put  $r(z) = 2 {\rm Re}(z_3) + |z_1^2 - z_2^3|^2$. Then
$\Delta_1^{\rm reg} (M,0) = 6$ but $\Delta_1(M,0) = \infty$. 
 \end{remark}

\section{proof of Theorem 1.1}

Following [D] or [D1],  we define both notions of order of contact.
Let ${\mathcal C}$ denote the collection of germs of nonconstant holomorphic maps 
$z:({\mathbb C}, 0) \to ( {\mathbb C}^n, 0)$. 
For a germ $r: ({\mathbb C}^n, 0) \to ({\mathbb R},0)$ of a smooth function we write $\nu(r)$ for its order of vanishing
at $0$. We also write $\nu(z)$ for the order of vanishing of a holomorphic map germ
$z \in {\mathcal C}$. Let ${\mathcal C}^*$ denote the collection of elements
$z$ in ${\mathcal C}$ for which $\nu(z)=1$.
As usual $z^*r$ denotes the (germ of a) map $t \mapsto r(z(t))$.

\begin{definition} Let $J$ be an ideal in $C^\infty_p$. We put
$\Delta_1(J) = \sup_{z \in {\mathcal C}} \inf_{h \in J} {\nu (z^* h) \over \nu(z)}$.
We put $\Delta_1^{\rm reg} (J) = \sup_{z \in {\mathcal C^*}} \inf_{h \in J} \nu (z^* h)$. \end{definition}

When $(M,0)$ is the germ of a real hypersurface, the ideal $J$ of germs vanishing on it is principal.
When $r$ is a generator of this ideal, 
the infima in Definition 2.1 are attained when $h= r$. When the supremum is finite,
only finitely many derivatives matter, and therefore many curves realize the supremum. 
When $(M,0)$ contains a unique holomorphic curve, 
the supremum is infinite and realized by only one curve. An example is given by
$r(z)= 2{\rm Re}(z_3) + |z_1^2 - z_2^3|^2$. The unique curve is given by $z(t)=(t^3,t^2, 0)$.

We write $\Delta_1(M,0)$ for $\Delta_1(J)$. We now prove Theorem 1.1.

\begin{proof} The proof combines property {\bf PS} with a Faa di Bruno formula
for powers of the Laplacian of the composite function $z^*r$.

Let $p$ be the origin in ${\mathbb C}^n$ and let $r$ generate the ideal of germs of smooth functions vanishing
on $(M,0)$. We may choose coordinates such that $\zeta=(z_1,...,z_{n-1})$ and 
$$ r(z) = 2 {\rm Re}(z_{n}) + g(\zeta) + 2 {\rm Im}(z_n) h(\zeta, {\rm Im}(z_n))) $$
where $g$ has no pure terms in its Taylor series up to as high an order as we desire.

There are many curves of maximum order of contact, that is, achieving the supremum in Definition 2.1.
It is noted in [D] that
one of these curves will satisfy $z_n(t) = 0$. The problem therefore reduces to showing the following. 
If there is a singular curve of multiplicity $m$ and with contact $4m$,
then there is a non-singular curve with contact $4$. Assume $g(\zeta)$ satisfies {\bf PS}. 
Let $z$ be a curve for which $\nu(z) = m$. Suppose that $\nu(z^*r) = 4m$. We will find a curve $\eta$
with $\nu(\eta) = 1$ and $\nu(z^*r)= 4$. 

We may assume that $({d \over dt})^k (z^*g)(0) = 0$ for $k \le 4m$.
Since {\bf PS} holds, and we are assuming that $\nu(z^*r) = 4m$, we know that
$L^{2m}(z^*g) (0) \ne 0$. We also are assuming that $z^{(j)}(0)=0$ for $0 \le j \le m-1$, but $z^{(m)}(0) \ne 0$.

In the following we will write $z^{(j)}$ to denote the $j$-th derivative of $z$ evaluated at $0$.
When $j=1$ we write $z'$ and when $j=2$ we write $z''$. Thus $z'$ and $z''$ are constant vectors.
Let $D^{ab}$ denote the symmetric multilinear form of type $(a,b)$ defined by
the derivatives of $g$ at the origin. Here there are $a$ holomorphic derivatives and $b$ barred derivatives.
Thus, for example,
$$ D^{11}(z',{\overline z}') = \sum_{j,k=1}^n g_{j {\overline k}} z'_j(0) {\overline {z'_k}}(0). $$
$$ D^{21} (z',z', {\overline z}') = \sum_{j,k,l=1}^n g_{jk {\overline l}} z_j'(0)z_k'(0) {\overline z}'_l (0). $$
We will work only with the notation on the left-hand side of these formulas.

We compute the first and second powers of the Laplacian on $z^*g$:
$$L(z^* g)(0) = D^{11}(z',{\overline z}') \eqno (3) $$
$$L^2(z^*g)(0) = D^{11}( z'', {\overline z}'') + D^{12}(z'', {\overline z}', {\overline z}') + D^{21}(z',z',{\overline z}'') + D^{22}( z',z', {\overline z}', {\overline z}'). \eqno (4) $$

When we apply $L$ an additional time, we begin to see combinatorial coefficients (related to Stirling numbers).
For the third power, we write only those terms involving $D^{11}, D^{12}, D^{21}, D^{22}$ because the others
will vanish in the crucial calculation below. Using the product and chain rules from calculus, we obtain
 
$$ L^3(z^*g) =  D^{11}( z^{(3)}, {\overline z}^{(3)}) + 3 D^{12}(z^{(3)}, {\overline z}^{(1)}, {\overline z}^{(2)}) + 
3 D^{21}(z^{(1)},z^{(2)} ,{\overline z}^{(3)})$$
$$ + 9 D^{22}( z^{(1)} ,z^{(2)}, {\overline z}^{(1)}, {\overline z}^{(2)}) + ..., \eqno (5) $$
where $...$ denotes terms involving $D^{ab}$ for $\max (a,b) \ge 3$. 

Let $z$ be a curve with $\nu(z) = m$ for $m\ge 2$. Then $z^{(j)} = 0$ for $0 \le j \le m-1$. 
Although $L^{2m}(z^*g)$ is an elaborate formula involving $D^{ab}$ for $1\le a,b \le 4m$,
the multi-linearity guarantees that any term with a slot $z^{(j)}$ for $j < m$ must vanish.
The total number of derivatives taken is $4m$. The only possible ways to get a weighted total of $4m$ derivatives
are those listed. We obtain
$$ L^{2m}(z^*g)(0) = $$
$$ D^{11}( z^{(2m)}, {\overline z}^{(2m)}) + 3 D^{12}(z^{(2m)}, {\overline z}^{(m)}, {\overline z}^{(m)}) + 3 D^{21}(z^{(m)},z^{(m)} ,{\overline z}^{(2m)})$$
$$ + 9 D^{22}( z^{(m)} ,z^{(m)}, {\overline z}^{(m)}, {\overline z}^{(m)}). \eqno (6)$$

Assume the expression in (6) is not zero.  Define $\eta$ by
$$ \eta(t) = {z^{(m)} \over \sqrt{3}}t + z^{(2m)} { t^2 \over 2} + ... .$$
Then $\eta' = {z^{(m)} \over \sqrt{3}}$ and $\eta'' = z^{(2m)}$. 
Then $\eta$ lies in ${\mathcal C}^*$. By (6), $\nu(\eta^*g) = 4$.
\end{proof}

\section{acknowledgements}

The author thanks the referee for suggesting some clarifications.
The author particularly thanks Jeff McNeal for noting the author's omission in [D]
and for sharing versions of the preprint [MM] with him.
The author ackowledges useful discussions with Dmitri Zaitsev, Siqi Fu, and Ming Xiao.
The important preprint [Z] by Zaitsev makes a systematic study
of fourth order invariants, but it does not include our Theorem 1.1.
The author acknowledges support from NSF Grant DMS 13-61001.

\section{bibliography}

\medskip

[BN]  V. Brinzanescu and A. Nicoara, On the relationship between D'Angelo q-type and Catlin q-type, {\it J. Geom. Anal.} 25 (2015), no. 3, 1701-1719. 

\medskip

[BN1] V. Brinzanescu, Vasile and A. Nicoara, Relating Catlin and D'Angelo q-types,
Math arXiv:1707.08294.

\medskip

[BS] H. P. Boas and E. J. Straube, On equality of line type and variety type of real hypersurfaces in 
${\mathbb C}^n$, {\it J. Geom. Anal.} 2 (1992), no. 2, 95-98. 

\medskip

[C] D. Catlin, Subelliptic estimates for the ${\overline \partial}$-Neumann problem on pseudoconvex domains,
{\it Ann. of Math.} (2) 126 (1987), no. 1, 131-191.

\medskip

[D] J. P. D'Angelo,  Several Complex Variables and the Geometry of Real Hypersurfaces,
CRC Press, Boca Raton, Fla., 1992.

\medskip

[D1] J. P.  D'Angelo, Real hypersurfaces, orders of contact, and applications, {\it  Annals of Math} (2) 115 (1982), no. 3, 615-637. 

\medskip

[DF]  K. Diederich and J. E. Fornaess, Pseudoconvex domains with real-analytic boundary, {\it Annals of Math} (2) 107 (1978), no. 2, 371-384. 

\medskip

[F] M. Fassina, The relationship between two notions of order of contact, preprint.

\medskip

[FIK] S. Fu, A. Isaev, and S. Krantz, 
Finite type conditions on Reinhardt domains, 
{\it Complex Variables Theory Appl.} 31 (1996), no. 4, 357-363.

\medskip

[K] J. J. Kohn, Subellipticity of the ${\overline \partial}$-Neumann problem on pseudoconvex domains: sufficient conditions, {\it  Acta Math} 142 (1979), no. 1-2, 79-122. 

\medskip

[M] J. D. McNeal, Convex domains of finite type, {\it J. Funct. Anal.} 108 (1992), no. 2, 361-373. 

\medskip

[MM] L. Mernik and J. McNeal, Regular versus singular order of contact on pseudoconvex hypersurfaces,
 Math arXiv:1708.02673.
\medskip 

[Z] D. Zaitsev, A geometric approach to Catlin's boundary systems, Math arXiv: 1704.01808.

\end{document}